\documentclass[12pt,a4paper]{article}
\usepackage{amsfonts}

\usepackage{amssymb}

\usepackage{amsmath}

\newtheorem{theorem}{Theorem}[section]

\newtheorem{lemma}[theorem]{Lemma}
\newtheorem{corollary}[theorem]{Corollary}

\date{}

\begin{document}

\title{Gr\"{o}bner-Shirshov bases for Rota-Baxter
algebras\footnote{Supported by the NNSF of China (No.10771077) and
the NSF of Guangdong province (No.06025062).} }

\author{
L. A. Bokut\footnote {Supported by RFBR 01-09-00157, LSS--344.2008.1
and SB RAS Integration grant No. 2009.97 (Russia).} \\
{\small \ School of Mathematical Sciences, South China Normal
University}\\
{\small Guangzhou 510631, P. R. China}\\
{\small Sobolev Institute of Mathematics, Russian Academy of
Sciences}\\
{\small Siberian Branch, Novosibirsk 630090, Russia}\\
{\small bokut@math.nsc.ru}\\
\\
 Yuqun
Chen\footnote {Corresponding author.} \ \ and Xueming Deng\\
{\small \ School of Mathematical Sciences, South China Normal
University}\\
{\small Guangzhou 510631, P. R. China}\\
{\small yqchen@scnu.edu.cn}\\
{\small dengxm860416@yhaoo.cn }} \vspace{4mm}

\maketitle

\noindent\textbf{Abstract:} In this paper, we establish the
Composition-Diamond lemma for associative nonunitary Rota-Baxter
algebras with weight $\lambda$. As applications, we obtain a linear
basis of a free commutative Rota-Baxter algebra without unity and
show that every countably generated Rota-Baxter algebra with weight
0 can be embedded into a two-generated Rota-Baxter algebra.



\section{Introduction}
Gr\"{o}bner bases and Gr\"{o}bner-Shirshov bases theories were
invented independently by A.I. Shirshov for non-associative algebras
and commutative (anti-commutative) non-associative algebras
\cite{Sh62a}, for Lie algebras (explicitly) and associative algebras
(implicitly) \cite{Sh62b}, for infinite series algebras (both formal
and convergent) by H. Hironaka \cite{Hi64} and for polynomial
algebras by B. Buchberger (first publication in \cite{Bu70}).
Gr\"{o}bner bases and Gr\"{o}bner-Shirshov bases theories have been
proved to be very useful in different branches of mathematics,
including commutative algebra and combinatorial algebra, see, for
example, the books
 \cite{AL, BWK93, BKu94, BuCL, BuW, CLO, Ei}, the papers \cite{Be78, Bo72,Bo76},
 and the surveys \cite{BC, BFKK00, BK03, BK05}.

The Shirshov's Composition-Diamond lemma--Buchberger's theorem is
the corner stone of the theories. This proposition says that in
appropriate free algebra $A_k(X)$ over a field $k$ with a free
generating set $X$ and a fixed monomial ordering, the following
conditions on a subset $S$ of $A_k(X)$ are equivalent: (i) Any
composition ($s$-polynomial) of  polynomials from $S$ is trivial;
(ii) If $f\in Id(S)$, then the maximal monomial $\bar{f}$ contains
some maximal monomial $\bar s$, where $s\in S$; (iii) The set
$Irr(S)$ of all (non-associative in general) words in $X $, which do
not contain any maximal word $\bar s, s\in S$, is a linear $k$-basis
of the algebra $A(X|S)=A(X)/Id(S)$ with generators $X$ and defining
relations $S$ (for Lie algebra case, $Irr(S)$ is a set of
Lyndon--Shirshov Lie words whose associative support do not contain
maximal associative words of polynomials from $S$).

Up to now, different versions of Composition-Diamond lemma are known
for the following classes of algebras apart those mentioned above:
(color) Lie super-algebras (\cite{Mik89, Mik92}) \cite{Mik96}, Lie
$p$-algebras \cite{Mik92}, associative conformal algebras
\cite{BFK}, modules \cite{KL}, dialgebras \cite{BCL08}, associative
algebras with multiple operators \cite{BCQ08}. In this paper, we
establish the Composition-Diamond lemma for associative nonunitary
Rota-Baxter algebras with weight $\lambda$. As applications, we
obtain a linear basis of a free commutative Rota-Baxter algebra
without unity and show that every countably generated Rota-Baxter
algebra with weight 0 can be embedded into a two-generated
Rota-Baxter algebra. In \cite{BFK}, a $\frac{1}{2}$-PBW Theorem for
associative conformal algebras was proved by L.A. Bokut, Y. Fong and
W.-F. Ke. Here we prove the similar results for dendriform dialgebra
and trialgebra as another two applications of the
Composition-Diamond lemma for associative nonunitary Rota-Baxter
algebras with weight $\lambda$.

Rota-Baxter algebras were invented by G. Baxter \cite{Bax} and
studied by G.-C. Rota \cite{R069, RO95, RO98} , F.V. Atkinson
\cite{At63} and P. Cartier \cite{Ca}. Since then, it has been
related to many topics in mathematics and mathematical physics, see,
for example, \cite{CK98, CK00, CK01, Der73, EF02, EGP, EG07, EG08c,
EGK04,  G02, G04, GS, GX, LHB,  Mi66, Mi69, N76}.

Not so many examples of Rota-Baxter algebras are known. They are
mostly free associative (commutative and with unity) Rota-Baxter
algebras \cite{AM, EG08a, EG08b} (\cite{GKg00a, GKg00b, GX}) and
low-dimensional Rota-Baxter algebras \cite{AB07, LB}. This paper
provides a systematical method to deal with Rota-Baxter algebras
defined by generators and defining relations.

In this paper, $k$ is a field of characteristic zero and Rota-Baxter
algebra always means Rota-Baxter algebra without unity.

We thank Ms Qiuhui Mo for some useful discussions.

\section{ Free Rota-Baxter algebra}

Let $A$ be an associative algebra over $k$ and $\lambda\in k$. Let a
$k$-linear operator $P: A\rightarrow A$ satisfy
$$
P(x)P(y)=P(P(x)y)+P(xP(y))+\lambda P(xy), \forall x,y\in A.
$$
Then $A$ is called a Rota-Baxter algebra with weight $\lambda$.

A subset $I$ of a Rota-Baxter algebra $A$ is called a Rota-Baxter
ideal of $A$ if $I$ is an ideal of $A$ such that $P(I)\subseteq I.$

Let $A, B$ be two Rota-Baxter algebras and $f: A\longrightarrow B$
be a mapping. Then $f$ is called a homomorphism if $f$ is a
$k$-algebra homomorphism such that for any $a\in A,
f(P(a))=P(f(a)).$

Recall that a free Rota-Baxter algebra with weight $\lambda$ on a
set $X$ is a Rota-Baxter algebra $A$ generated by $X$ with a natural
mapping $i: X\rightarrow A$ such that, for any Rota-Baxter algebra
$R$ with weight $\lambda$ and any map $f: X\rightarrow R$, there
exists a unique homomorphism $\tilde{f}: A\rightarrow R$ such that
$\tilde{f}\cdot i=f$.

The free Rota-Baxter algebra generated by a nonempty set $X$ is
given by K. Ebrahimi-Fard and L. Guo \cite{EG08a}.

Let $X$ be a nonempty set, $S(X)$ the free semigroup generated by
$X$ without identity and $P$ a symbol of a unary operation. For any
two nonempty sets $Y$ and $Z$, denote by
$$
{\Lambda}_{P}{(Y,Z)}=(\cup_{r\geq 0}(YP(Z))^{r}Y)\cup(\cup_{r\geq
1}(YP(Z))^{r})\cup(\cup_{r\geq 0}(P(Z)Y)^{r}P(Z))\cup(\cup_{r\geq
1}(P(Z)Y)^{r}),
$$
where for a set $T$,  $T^0$ means the empty set.

 \noindent{\bf Remark:} In
${\Lambda}_{P}{(Y,Z)}$, there are no words with a subword
$P(z_1)P(z_2)$ where $z_1,z_2\in Z.$

Define
\begin{eqnarray*}
\Phi_{0}&=&S(X)\\
\vdots\ \ & &\ \ \ \ \vdots\\
\Phi_{n}&=&{\Lambda}_{P}(\Phi_{0},\Phi_{n-1}) \\
\vdots\ \ & &\ \ \ \ \vdots\\
\end{eqnarray*}
Then
$$
\Phi_{0} \subset \cdots \subset\Phi_{n} \subset
 \cdots
$$
Let
$$
\Phi(X)=\cup_{n\geq 0} \Phi_{n}.
$$
Clearly, $P(\Phi(X))\subset \Phi(X)$. If $u\in X\cup P({\Phi}(X))$,
then $u$ is called prime.  For any $u\in \Phi(X)$, $u$ has a unique
form
 $u=u_{1}u_{2}\cdots u_{n}$ where $u_{i}$  is prime, $i=1, 2,\dots,n$, and
 $u_i,u_{i+1}$ can not both have forms as $p(u'_i)$ and
 $p(u'_{i+1})$. If this is the case, then we define the breath of $u$ to be $n$, denoted by
 $bre(u)=n$.

For any $u\in \Phi(X)$ and for a set $T\subseteq X\cup \{P\}$,
denote by $deg_{T}(u)$ the number of occurrences of $t\in T$ in $u$.
 Let
$$
Deg(u)=(deg_{\{P\}\cup X}(u), deg_{\{P\}}(u)).
$$
We order $Deg(u)$
lexicographically.

Let $k\Phi(X)$ be a free  $k$-module with $k$-basis $\Phi(X)$ and
$\lambda \in k$ a fixed element. Extend linearly $P: \
k\Phi(X)\rightarrow k\Phi(X),\ u\mapsto P(u)$ where $u\in \Phi(X)$.

Now we define the multiplication in $k\Phi(X)$.

Firstly, for $u,v\in X\cup P(\Phi(X))$, define

$u\cdot v= \left\{ \begin{array}
              {l@{\quad}l}
              P(P(u')\cdot v')+P(u'\cdot P(v'))+\lambda P(u'\cdot v'), & \mbox{if}\ u=P(u'),v=P(v'); \\
              uv, & \mbox{otherwise}.
              \end{array} \right. $

Secondly, for any $u=u_{1}u_{2}\cdots u_{s}, v=v_{1}v_{2}\cdots
v_{l}\in \Phi(X)$ where $u_{i}, v_{j}$ are prime, $i=1,2,\dots,s,
j=1,2,\dots,l$, define
$$
u\cdot v=u_{1}u_{2}\cdots u_{s-1}(u_{s}\cdot v_{1})v_{2}\cdots
v_{l}.
$$

Equipping with the above concepts, $k\Phi(X)$ is the free
Rota-Baxter algebra with weight $\lambda$ generated by $X$, see
\cite{EG08a}.

We denote by $RB(X)$ the free Rota-Baxter algebra with weight
$\lambda$ generated by $X$.

\section{Composition-Diamond lemma for Rota-Baxter algebras}

In this section, we establish the Composition-Diamond lemma for
Rota-Baxter algebras with weight $\lambda$.

Let $\cal N$ be the set of  non-negative integers and  $\cal N^+$
the set of positive integers.

Let the notations be as before. We have to order $\Phi(X)$. Let $X$
be a well ordered set. Let us define an ordering $>$ on $\Phi(X)$ by
induction on the $Deg$-function.

For any $u,v\in \Phi(X)$, if $Deg(u)>Deg(v)$, then $u>v$.

If $Deg(u)=Deg(v)=(n,m)$, then we define $u>v$ by induction on
$(n,m)$.

If $(n,m)=(1,0),$ then $u,v\in X$ and we use the ordering on $X$.
Suppose that for $(n,m)$ the ordering is defined where
$(n,m)\geq(1,0)$. Let $(n,m)<(n',m')=Deg(u)= Deg(v)$.  If $u,v \in
P(\Phi(X))$, say $u=P(u')$ and $v=P(v')$, then $u>v$ if and only if
$u'>v'$ by induction. Otherwise $u=u_1u_2\cdots u_l$ and
$v=v_1v_2\cdots v_s$ where $l>1$ or $s>1$, then $u>v$ if and only if
$(u_1, u_2, \dots, u_l)>(v_1, v_2, \dots, v_s )$ lexicographically
by induction.

It is clear that $>$ is a well ordering on $\Phi(X)$. Throughout
this paper, we will use this ordering.

Now, for any $0\neq f\in RB(X)$, $f$ has the leading term $\bar{f}$
and $f=\alpha_{1}\bar{f}+\sum_{i=2}^{n}\alpha_iu_i$ where
$\bar{f},u_i\in \Phi(X), \bar{f}>u_i, 0\neq\alpha_1, \alpha_i \in
k$. Denote by $lc(f)$ the coefficient of the leading term $\bar{f}$.
If $lc(f)=1$, we call $f$ monic.

For any $u\in \Phi(X)$ and any non-negative integer $n$, we denote
by
$$
P^{n}(u)=\underbrace{P(P(\cdots P}_n(u)\cdots)) \mbox{ and }
P^0(u)=u.
$$

Clearly, if $u\in \Phi(X)$, then there exists a unique $n\in\cal N$
such that $u=P^n(u')$ where $u'\in \Phi(X)\backslash P(\Phi(X))$.

The proof of the following lemma is straightforward. We omit
details.

\begin{lemma}\label{Deng1}
For any $u,v\in \Phi(X),\  n, m\geq 1$, we have
\begin{eqnarray*}
&&P^{n}(u)\cdot P^{m}(v)\\
&=&\sum_{s=1}^{n}\alpha_{(n,m,s)}P^{n+m-s}(P^{s}(u)\cdot v)
+\sum_{l=1}^{m}\beta_{(n,m,l)}P^{n+m-l}(u\cdot P^{l}(v))+\lambda
{\varepsilon}(n,m,u,v)
\end{eqnarray*}
where $\alpha_{(n,m,s)},\  \beta_{(n,m,l)} \in {\cal N^+}, \
{\varepsilon}(n,m,u,v)\in P(RB(X)),
deg_{\{P\}}(\overline{{\varepsilon}(n,m,u,v)})=
deg_{\{P\}}(P^n(u))+deg_{\{P\}}(P^m(v))-1\ \mbox{and}\
deg_{X}(\overline{{\varepsilon}(n,m,u,v)})=deg_{X}(u)+deg_{X}(v) .$
Moreover, if we define $\alpha_{(n,0,s)}=0, \beta_{(0,m,l)}=0$, then
the coefficients satisfy the recursive relations:
\begin{enumerate}
\item[1)]\ $\alpha_{(n,1,s)}=1\ \mbox{and}\ \beta_{(1,m,l)}=1\ \mbox{where}\ 1\leq s\leq n-1, 1\leq l\leq m-1.$
\item[2)]\ $\alpha_{(n,m,s)}=\alpha_{(n-1,m,s)}+\alpha_{(n,m-1,s)}\
\mbox{and}\ \beta_{(n,m,l)} =\beta_{(n-1,m,l)} +\beta_{(n,m-1,l)}\
\mbox{where}\ 1\leq s\leq n-1, 1\leq l\leq m-1.$
\item[3)]\ $\alpha_{(n,m,n)}=\beta_{(n,m,m)}=1.$
\end{enumerate}
\end{lemma}

\begin{lemma}\label{Deng2}
For any $u,v\in \Phi(X)$, if $u,v\notin P(\Phi(X))$, then
$\overline{P(u)\cdot v}> \overline{u\cdot P(v)}$. Moreover, for any
$m\geq1, n\geq1, \ \overline{P^{n}(u)\cdot
P^{m}(v)}=P^{n+m-1}({\overline{P(u)\cdot v}})$.
\end{lemma}
{\bf Proof:} There are four cases to consider.

Case 1: $bre(u)=bre(v)=1$, i.e., $u=x, v=y \in X$.

Case 2: $bre(u)=1$ and $bre(v)=l>1$, i.e.,  $u=x\in  X,
v=v_{1}v_{2}\cdots v_{l}$.

Case 3:  $bre(u)=s>1$  and $bre(v)=1$, i.e., $u=u_{1}\cdots
u_{s-1}u_{s},\ v=y\in X$.

Since Case 1, Case 2 and Case 3 are simple to prove, we just prove
the following.

Case 4: $bre(u)=s>1$ and $bre(v)=l>1$, say $u=u_{1}\cdots
u_{s-1}u_{s}, v=v_{1}v_{2}\cdots v_{l}$. Then we have the following
two subcases to consider.
\begin{description}
\item{(i)} If $v_{1}=y\in X$, then $P(u)\cdot v=P(u_{1}\cdots
u_{s-1}u_{s})yv_{2}\cdots v_{l}, u\cdot P(v)=u_{1}\cdots
u_{s-1}u_{s}\cdot P(yv_{2}\cdots v_{l})$. Since $P(u_{1}\cdots
u_{s-1}u_{s})>u_{1}$,
$$
\overline{P(u)\cdot v}=P(u_{1}\cdots u_{s-1}u_{s})y v_{2}\cdots
v_{l}>u_{1}\cdots u_{s-1} \overline{u_{s}\cdot P(yv_{2}\cdots
v_{l})}=\overline{u\cdot P(v)}.
$$
\item{(ii)} If $v_{1}=P(v_{1}')$, then $P(u)\cdot v=(P(u)\cdot
P(v_{1}'))v_{2}\cdots v_{l}$ and  $u\cdot P(v)=u_{1}\cdots
u_{s-1}(u_{s}\cdot P(v)).$ Since $\overline{P(u)\cdot P(v_{1}')}\in
P(\Phi(X))$ and $\overline{P(u)\cdot P(v_{1}')}>u_{1}$,
$$
\overline{P(u)\cdot v}=\overline{P(u)\cdot P(v_1')}v_{2}\cdots
v_{l}>u_{1}\cdots u_{s-1}\overline{u_{s}\cdot P(P(v_{1}')v_{2}\cdots
v_{l})}=\overline{u\cdot P(v)}. \ \ \ \ \ \ \ \hfill \blacksquare
$$
\end{description}

\begin{lemma}\label{Deng3}
For any $u,v\in \Phi(X)$, if $u>v$,  then $\overline{u\cdot w}>
\overline{v\cdot w}\ \  and \ \ \overline{w\cdot u}>
\overline{w\cdot v}$ for any $w\in \Phi(X).$
\end{lemma}

{\bf Proof:}  Firstly, we show $\overline{u\cdot w}>
\overline{v\cdot w}$. We may assume $u=u_{1}\cdots u_{s},
v=v_{1}\cdots
 v_{l}, s,l\geq 1$ and $u_{1}>v_{1}$ where $u_i,v_j$
are prime, $i=1,2,\dots,s, j=1,2,\dots,l$, and $Deg(u)=Deg(v)$. We
may also assume that $w=P(w')\in P(\Phi(X))$.

It is noted that if $s>1$ and $l=1$, then, since $Deg(u)=Deg(v)$, we
will have $v>u$, a contradiction. So we just need to consider the
following cases.

Case 1: $s>1,\ l>1$ and Case 2: $s=1,\ l>1$. These two cases are
clear.

Case 3: $s=l=1$. If $u_{1}\in X$, then $v_{1}\in X$ and the result
holds clearly. If $u_{1}\in P(\Phi(X))$, then $v_{1}\in P(\Phi(X))$.
We may assume that $u_{1}=P^{n}(u_{1}'), v_{1}=P^{n'}(v_{1}')$ with
$u_{1}', v_{1}'\notin P(\Phi(X))$ and $w=P^{m}(w''), w''\notin
P(\Phi(X))$. Then we need to consider two subcases.
\begin{description}
\item (i) If $n>n'$, then $\overline{u\cdot
w}=P^{n+m-1}(\overline{P(u_{1}')\cdot
w''})>P^{n'+m-1}(\overline{P(v_{1}')\cdot w''})=\overline{v\cdot w}$
by Lemma \ref{Deng2}.
\item (ii) If $n=n'$, then $u_1'>v_1'$.
By Lemma \ref{Deng2} we have
$$
\overline{u\cdot w}=P^{n+m-1}(\overline{P(u_{1}')\cdot w''}),\ \
\overline{v\cdot w}=P^{n+m-1}(\overline{P(v_{1}')\cdot w''}).
$$
By induction on $deg_{\{P\}\cup X}(w)$, $\overline{P(u_{1}')\cdot
w''}>\overline{P(v_{1}')\cdot w''}$ and so $\overline{u\cdot
w}>\overline{v\cdot w}$.
\end{description}

This shows that $\overline{u\cdot w}>\overline{v\cdot w}.$

Secondly, we show that $\overline{w\cdot u}> \overline{w\cdot v}$.

We may assume that $w=P(w')\in P(\Phi(X))$, \ $Deg(u)=Deg(v),\
u=u_{1}u_{2}\cdots u_{s}> v=v_{1}v_{2}\cdots v_{l}$ where $ u_{i},
v_{j}$ are prime, $i=1,2,\dots,s, j=1,2,\dots,l$.  Then we have the
following two cases.

Case 1: $u_{1}>v_{1}$. We only need to consider the case of
$u_{1}=P^{n}(u_{1}')\in P(\Phi(X)),\ v_{1}=P^{n'}(v_{1}' )\in
P(\Phi(X))$ where $n,n'\geq 1, u_{1}', v_{1}'\notin P(\Phi(X)).$ Let
$w=P^{m}(w'')\in P(\Phi(X))$ where $m\geq 1\ \mbox{and}\ w''\notin
P(\Phi(X))$. Then we have the following subcases to consider.
\begin{description}
\item{(i)} If $n>n'$, then $\overline{w\cdot
u}=P^{n+m-1}(\overline{P(w'')\cdot u_{1}'})u_{2}\cdots
u_{s}>P^{n'+m-1}(\overline{P(w'')\cdot v_{1}'}) v_{2}\cdots
v_{l}=\overline{w\cdot v}$ by Lemma \ref{Deng2}.
\item{(ii)} If $n=n'$, then $u_{1}'>v_{1}'.$ By induction on
$deg_{\{P\}\cup X}(u)+deg_{\{P\}\cup X}(v)$, we have
$\overline{P(w'')\cdot u_{1}'}>\overline{P(w'')\cdot v_{1}'}$.
Therefore, by Lemma \ref{Deng2} $$\overline{w\cdot
u}=P^{n+m-1}(\overline{P(w'')\cdot u_{1}')}u_{2}\cdots
u_{s}>P^{n+m-1}(\overline{P(w'')\cdot v_{1}'}) v_{2}\cdots
v_{l}=\overline{w\cdot v}.$$
\end{description}
Case 2: $u_{1}=v_{1}$. This case is clear.  $ \hfill \blacksquare $

\ \

For any $u, v\in \Phi(X)$, clearly,
$$
u>v \Longrightarrow P(u)>P(v).
$$

Let $\star$ be a symbol and $\star\notin X$. By a
$\star$-Rota-Baxter word we mean any expression in $\Phi(X\cup
\{\star\})$ with only one occurrence of $\star$. The set of all
$\star$-Rota-Baxter words on $X$ is denoted by $\Phi^\star (X)$.

Let $u$ be a $\star$-Rota-Baxter word and $s\in RB(X)$. Then we call
$$
u|_{s}=u|_{\star\mapsto s}
$$
an $s$-Rota-Baxter word. For short, we  call $u|_{s}$ an $s$-word.

In other words, an $s$-word $u|_s$ means that we have replaced the
$\star$ of $u$ by $s$.

For example, if $u=P(x_1)x_2P^2(\star)x_4P(x_5)$, then
$u|_s=P(x_1)x_2P^2(s)x_4P(x_5)$ is an $s$-word.

Similarly, we can define  $(\star_1, \star_2)$-Rota-Baxter words as
expressions in  $\Phi(X\cup \{\star_1, \star_2\})$  with only one
occurrence of $\star_1$ and only one occurrence of $\star_2$. Let us
denote by $\Phi^{\star_1, \star_2} (X)$ the set of all $(\star_1,
\star_2)$-Rota-Baxter words. Let $u\in \Phi^{\star_1, \star_2} (X)$.
Then we call
$$
u|_{s_1, s_2}= u|_{\star_1\mapsto s_1,\star_2\mapsto s_2}
$$
an $s_1$-$s_2$-Rota-Baxter word. For short, we call $u|_{s_1, s_2}$
an $s_1$-$s_2$-word.

If $\overline{u|_{s}}=u|_{\overline s}$, then we call $u|_{s}$ a
normal $s$-word.

By Lemma \ref{Deng3}, we have the following lemma which shows that
the ordering $>$ on $\Phi(X)$ is monomial.

\begin{lemma}\label{Deng4}
For any $u,v\in \Phi(X),\ w\in\Phi^{\star}(X)$,
$$
u>v\Longrightarrow \overline{w|_{u}}>\overline{w|_{v}}
$$ where $w|_u=w|_{\star\mapsto u}$ and $w|_v=w|_{\star\mapsto v}$.$\hfill
\blacksquare $
\end{lemma}

\noindent{\bf Remark:} If $u|_{s}$ is a normal $s$-word, then
$P^{l}(u|_{s})$ is also a normal $s$-word where $l\in {\cal N}.$

\ \

In order to describe the ideal $Id(S)$ of $RB(X)$ generated by $S$,
we introduce the concept of $P$-$s$-words.

For any nonempty set $X$, define
\begin{eqnarray*}
\Psi_{0}&=&S(X)\\
\Psi_{1}&=&S(X\cup P(S(X))\\
\vdots\ \ & &\ \ \ \ \vdots\\
\Psi_{n}&=&S(X\cup P(\Psi_{n-1})) \\
\vdots\ \ & &\ \ \ \ \vdots
\end{eqnarray*}
Then we have
$$
\Psi_{0} \subset \Psi_{1} \subset \cdots \subset\Psi_{n} \subset
 \cdots
$$
Let
$$
\Psi(X)=\cup_{n\geq 0} \Psi_{n}.
$$

Clearly, $\Phi(X)\subset\Psi(X)$.

\ \

 Now, we can define
$\Psi^{\star}(X)$ as the definition of $\Phi^{\star}(X)$.

Let $S\subseteq RB(X)$. Then it is clear that  $Id(S)=span_k\{w|_s\
|\ w\in \Psi^{\star}(X), s\in S \}$ where $span_k\{w|_s\ |\ w\in
\Psi^{\star}(X), s\in S \}$ is the $k$-subspace of $RB(X)$ generated
by $\{w|_s\ |\ w\in \Psi^{\star}(X), s\in S \}$.

Let $ w\in \Psi^{\star}(X)$ and $s\in S$. Then we call $w|_s$ a
$P$-$s$-word.

The following lemma says that each $P$-$s$-word is a linear
combination of $s$-words.

\begin{lemma}\label{lemma3.5}
Let $S\subseteq RB(X),\ u'\in \Psi^{\star}(X)$ and $s\in S$. Then
the $P$-$s$-word $u'|_s$ has an expression:
$$
u'|_s=\sum\alpha_iu_i|_s, \ \ \ \mbox{where each} \ \ \ \alpha_i\in
k,\ u_i|_s \mbox{ is }\ s\mbox{-word and }\
\overline{u_i|_s}\leq\overline{u'|_s}.
$$
In particular, $Id(S)=span_k\{w|_s\ |\ w\in \Phi^{\star}(X), s\in S
\}.$
\end{lemma}
{\bf Proof:}
 For any $u\in \Psi(X)$, in $\Psi(X)$, $u$ has a unique
expression
$$
u=u_1\cdots u_n, \ \ \ \ u_i\in X\cup P(\Psi(X)).
$$
Such an $n$ is called the $P$-length of $u$.

Let $u'|_s$ be a $P$-$s$-word. Denote by $u''=u'|_{\star\mapsto
x}=u_1\cdots u_n$, where $x\in X,\ u_i\in X\cup P(\Psi(X))$. We
prove the result by induction on the $P$-length $n$ of $u''$ and on
$\overline{u'|_s}$.

We may assume that $u''\not\in \Phi(X)$.

Case 1. $n=1$. Then $u'|_s=P(u|_s)$. Since
$\overline{u|_s}<\overline{u'|_s}$, the result follows by induction
on $\overline{u'|_s}$.

Case 2. $n\geq 2$. Then there are two cases to consider.
\begin{description}

\item{}2.1 There exists $1\leq i\leq n-1$ such that $u_i=P(u'_i)$ and
$u_{i+1}=P(u'_{i+1})$ where $u'_i, u'_{i+1}\in \Psi(X)$. Thus, by
the Rota-Baxter identity, in $RB(X)$,
$$
u''=u_1\cdots u_{i-1}(P(P(u'_i)u'_{i+1})+P(u'_iP(u'_{i+1}))+\lambda
P(u'_iu'_{i+1}))u_{i+2}\cdots u_{n},
$$
where in the right side, each term has $P$-length $n-1$ and
$\overline{(u_1\cdots u_{i-1}P(P(u'_i)u'_{i+1})u_{i+2}\cdots
u_{n})|_s},\\ \overline{(u_1\cdots
u_{i-1}P(u'_iP(u'_{i+1}))u_{i+2}\cdots u_{n})|_s},\
 \overline{(u_1\cdots u_{i-1} P(u'_iu'_{i+1})u_{i+2}\cdots
u_{n})|_s}\leq \overline{u'|_s}$. Now the result follows by
induction on $n$.

\item{}2.2  For each $j,\ 1\leq j\leq n-1,\ \{u_j,u_{j+1}\}\cap
X\neq\emptyset$. Let $u'|_{s}=u_1\cdots u_{m-1}u_m|_su_{m+1}\cdots
u_n$.

For $u_m$, there are two cases to consider: $u_m=P(v)$ and $u_m=x$.

For each $u_i,\  1\leq i\leq n,\ i\neq m$, by the Rota-Baxter
identity, we may assume
$$
u_i=\sum \alpha_{iq_i}v_{iq_i},\ \ \ \ \  \mbox{ where each }
\alpha_{iq_i}\in k,\ v_{iq_i}\in \Phi(X),\
v_{iq_i}\leq\overline{u_i}.
$$
Thus,
$$
u_1\cdots u_m\cdots
u_n=\sum\alpha_{1q_1}\cdots\widehat{\alpha_{mq_m}}\cdots\alpha_{nq_n}v_{1q_1}\cdots
u_m\cdots v_{nq_n}
$$
where each $v_{1q_1} \cdots y\cdots v_{nq_n}\in \Phi^{\star}(X), \
y\in \{P(\star), \star\}$, $v_{1q_1}\cdots \overline{u_m}\cdots
v_{nq_n}\leq \overline{u_1\cdots u_m\cdots u_n}$, and $v_{1q_1}
\cdots{u_m} \cdots v_{nq_n}$ has  $P$-length $n$.

If $u_m=P(v)$ then  the result follows from the induction on
$\overline{u'|_s}$ since $\overline{v|_s}<\overline{u'|_s}$.

If $u_m=x$ then the result is clear.
\end{description}

The proof is complete.
  $\hfill \blacksquare$

 \ \

Let $f,g \in RB(X)$ be monic with $\overline{f}=u_{1}u_{2}\cdots
u_{n}$ where each $u_{i} $ is prime. Then, there are four kinds of
compositions.
\begin{description}
\item{(i)} If $u_{n}\in P(\Phi(X))$, then we define composition of
right multiplication  as\ \ $f\cdot u$ \ \ where $u\in P(\Phi(X)).$

\item{(ii)} If $u_{1}\in P(\Phi(X))$, then we define composition of
left multiplication as \ \ $u\cdot f$ \ \  where \ $u\in
P(\Phi(X)).$

\item{(iii)} If there exits a $w=\overline{f}a=b
\overline{g}$ \ where $fa$ is normal $f$-word and $bg$ is normal
$g$-word, $a, b \in \Phi(X)$ and $deg_{\{P\}\cup
X}(w)<deg_{\{P\}\cup X}(\overline{f})+deg_{\{P\}\cup
X}(\overline{g})$, then we define the intersection composition of
$f$ and $g$ with respect to $w$ as $(f,g)_{w}=f\cdot a-b\cdot g.$

\item{(iv)} If there exists a $w=\overline{f}=u|_{\overline{g}}$ where
$u\in \Phi^{\star}(X)$, then we define the inclusion composition of
$f$ and $g$ with respect to $w$ as $(f,g)_{w}=f-u|_{g}.$ Note that
if this is the case, then  $u|_g$ is a normal $g$-word.
\end{description}

We call $w$ in $(f,g)_{w}$ the ambiguity with respect to $f$ and
$g$. By Lemma \ref{Deng4},
$$
\overline{(f,g)_w}<w.
$$

Let $S\subset RB(X)$ be a set of monic polynomials. Then the
composition $(f,g)_w$ is called trivial modulo $(S,w)$ if
$$
(f,g)_w=\sum_{i}\alpha_iu_i|_{s_{i}}
$$
where each $\alpha_i\in k$,  $s_i\in S$, $u_i|_{s_i}$ is normal
$s_i$-word and $u_i|_{\overline{s_i}}< w$. If this is the case, then
we write
$$
(f,g)_w\equiv 0 \ \ mod (S,w).
$$

The composition of left (right) multiplication is called trivial
$mod(S)$ if
$$
u\cdot f=\sum_{i} \alpha_{i}u_i|_{s_{i}}\ \ (f\cdot u=\sum_{i}
\alpha_{i}u_i|_{s_{i}})
$$ where each $\alpha_i\in k$,  $s_i\in S$, $u_i|_{s_i}$ is normal $s_i$-word and
$u_i|_{\overline{s_i}}\leq \overline{u\cdot f}\ \
(u_i|_{\overline{s_i}}\leq \overline{f\cdot u})$. If this is the
case, then we write
$$
u\cdot f \equiv 0 \ \ mod(S)\ \ (f\cdot u \equiv 0 \ \ mod(S)).
$$

In general, for any two polynomials $p$ and $q$, $ p\equiv q \ \ mod
(S,w) $ means that $ p-q=\sum_{i}\alpha_iu_i|_{s_i}$ where each
$\alpha_i\in k$,  $s_i\in S$, $u_i|_{s_i}$ is normal $s_i$-word and
$u_i|_{\overline{s_i}}< w$.

$S$ is called a Gr\"{o}bner-Shirshov basis  in  $RB(X)$ if each
composition $(f,g)_w$ of $f,g\in S$  is trivial  $mod(S,w)$ and each
composition of left multiplication and right multiplication  is
trivial $mod(S)$.

\begin{lemma}\label{Deng5}
Let $S\subset RB(X)$ be a set of monic polynomials. If each
composition of left multiplication and right multiplication of $S$
is trivial $mod(S)$,  then each $s$-word $u|_{s}$ has an expression
$$
u|_{s}=\sum_{i} \alpha_{i}u_{i}|_{s_{i}}
$$
where each $\alpha_{i}\in k, s_{i}\in S, u_{i}|_{s_{i}}$ is normal
$s_{i}$-word and $u_{i}|_{\overline{s_{i}}}\leq \overline{u|_{s}}$.
\end{lemma}
{\bf Proof:} We prove the result by induction on $\overline{u|_s}$.
If $\overline{u|_s}=\overline{s}$, then there is nothing to prove.
Suppose that $\overline{u|_s}>\overline{s}$. We consider the
following two cases.

Case 1: $u|_{s}=w_{1}sw_{2}$ where $w_1, w_2 \in \Phi(X)$ ($w_{1}, \
w_{2}$ may be empty). We only consider the case $w_1, w_2 \in
\Phi(X)$. The other cases are simple.  By induction we assume that
$w_1s=\sum_i\alpha_{i}v_i|_{s_i}$ where each $v_i|_{s_i}$ is a
normal $s_i$-word, $s_i\in S, \alpha_{i}\in k$ and
$v_i|_{\overline{s_i}}\leq \overline{w_1s}$.

Suppose that $w_2=w_{21}w_{22}\cdots w_{2n}$ and
$\overline{s_i}=q_{i1}q_{i2}\cdots q_{it}$ where   $w_{2j}, \
q_{ij}$ are prime.

There are two cases to consider.

(I) If $v_i|_{s_i}=w_{1i}s_iw_{2i}$, then we need to consider the
following four subcases.

Subcase 1: $w_{1i}$ and $w_{2i}$ are both empty words. Then
$$
v_i|_{s_i}w_2=s_iw_2.
$$
 In this case, we also have the
following two subcases to consider.
\begin{description}
\item {(I-1-1)}:  $q_{it}\in X$. In this case, $v_i|_{s_i}w_2=s_iw_{2}$ is clearly a
normal $s_i$-word.
\item {(I-1-2)}:  $q_{it}\in P(\Phi(X))$.  If $w_{21}\in X$, then
$v_i|_{s_i}w_2=s_iw_{2}$ is normal $s_i$-word clearly. Otherwise,
since each composition of right multiplication is trivial,
$$
s_iw_{2}=\sum_{j}\gamma_{ij}v'_{ij}|_{s_{ij}}\ \ \ \ \ \ \ \ \ \ \ \
\ \ \ \ \ \ \  \ \ \ \ \ \ \ \ \ \ \ \ (*)
$$
where each
$\gamma_{ij}\in k, \ s_{ij}\in S,\ v'_{ij}|_{s_{ij}}$ is normal
$s_{ij}$-word and ${v'_{ij}|_{\overline{s_{ij}}}}\leq
\overline{s_iw_{2}}=\overline{v_i|_{s_i}w_2}\leq
\overline{w_1sw_2}=\overline{u|_{s}}$.
\end{description}

Subcase 2: $w_{1i}\in \Phi(X)$ and $w_{2i}$ is empty word. Then
$$
v_i|_{s_i}w_2=w_{1i}s_iw_2.
$$
\begin{description}
\item{} If $bre(\overline{s_i})=1$, then we have the following two subcases
to consider.
\begin{description}
\item{(I-2-1)}:  $\overline{s_i}=x\in X$. Then $v_i|_{s_i}w_2$
is a normal $s_i$-word clearly.
\item{(I-2-2)}:  $\overline{s_i}\in P(\Phi(X))$.  Then
we need to consider the following two subcases.
\begin{description}
\item{(i)}: If $w_{21}\in X$, then $v_i|_{s_i}w_2$
is a normal $s_i$-word.
\item{(ii)}: If $w_{21}\in P(\Phi(X))$, then because each composition of
right multiplication is trivial, we have (*). Since
$v_i|_{s_i}=w_{1i}s_i$ is normal $s_i$-word and $\overline{s_i}\in
P(\Phi(X))$, we have $w_{1i}=w_{1i}'x$ where $x\in X$. Therefore,
each $w_{1i}v'_{ij}|_{s_{ij}}$ is still normal $s_{ij}$-word. Thus,
$v_i|_{s_i}w_2=\sum_j\beta_{ij}w_{1i}v'_{ij}|_{s_{ij}}$ where each
$u_{ij}|_{s_{ij}}=w_{1i}v'_{ij}|_{s_{ij}}$ is normal $s_{ij}$-word
and $u_{ij}|_{\overline{s_{ij}}}=
w_{1i}v'_{ij}|_{\overline{s_{ij}}}\leq
\overline{w_{1i}\overline{s_iw_{2}}}=\overline{v_i|_{s_i}w_2}\leq\overline{u|_{s}}$.
\end{description}
\end{description}
\item If $bre(\overline{s_i})=t>1$, we also have the
following two subcases to consider.
\begin{description}
\item {(I-2-3)}:  $q_{it}\in X$. In this case, $v_i|_{s_i}w_2=w_{1i}s_iw_2$ is clearly a
normal $s_i$-word.
\item {(I-2-4)}:  $q_{it}\in P(\Phi(X))$.  If $w_{21}\in X$, then
$v_i|_{s_i}w_2=w_{1i}s_iw_2$ is normal $s_i$-word clearly.
Otherwise, since each composition of right multiplication is
trivial, we have (*). We need to consider the following two cases.
\begin{description}
\item{(i)}:  $q_{i1}=P(q_{i1}')$. Similar to the case (ii) in (I-2-2), we have the result.
\item{(ii)}:  $q_{i1}\in X$. We need to consider the following two
subcases.
\begin{description}
\item{(a)}: If $v'_{ij}|_{\overline{s_{ij}}}<\overline{s_iw_{2}}$, then
$\overline{w_{1i}v'_{ij}|_{s_{ij}}}<\overline{v_i|_{s_i}w_2}\leq\overline{u|_s}$.
By induction, $w_{1i}v'_{ij}|_{s_{ij}}=\sum_{m}
\alpha_{ijm}u_{ijm}|_{s_{ijm}}$ where each $s_{ijl}\in S,
u_{ijm}|_{s_{ijm}}$ is normal $s_{ijm}$-word and
$\overline{u_{ijm}|_{s_{ijm}}}\leq
\overline{w_{1i}v'_{ij}|_{s_{ij}}}<\overline{u|_{s}}$.
\item{(b)}:  If $v'_{ij}|_{\overline{s_{ij}}}=\overline{s_iw_{2}}$,
because $\overline{s_iw_{2}}=q_{i1}\cdots
q_{i(t-1)}\overline{q_{it}w_{2}}$ and $q_{i1}\in X$, we have
$w_{1i}v'_{ij}|_{s_{ij}}$ is still normal $s_{ij}$-word.
\end{description}
\end{description}
\end{description}
\end{description}

Subcase 3: $w_{1i}, w_{2i}\in \Phi(X)$ and Subcase 4: $w_{1i}$ is a
empty word and $w_{2i}\in \Phi(X)$. These two subcases are simple
to prove.

(II) $v_i|_{s_i}=w_{1i}P(v'_i|_{s_i})w_{2i}$. Since $v_i|_{s_i}$ is
a normal $s_i$ word, we may just consider the case $w_{1i}$ and
$w_{2i}$ are both empty word (other cases are simple). Moreover, we
may suppose $w_2=P(w_2')$. Then we have
$$
v_i|_{s_i}w_2=P(v'_i|_{s_i})\cdot P(w_2')=P(P(v'_i|_{s_i})\cdot
w_2')+P(v'_i|_{s_i}\cdot P(w_2'))+\lambda P(v'_i|_{s_i}\cdot w_2').
$$
Similar to the proof of Lemma \ref{lemma3.5}, we have the result.

This proves the Case 1.

Case 2: $u|_{s}=w_{1}P(u'|_{s})w_{2}$ where $ u'\in \Phi^{\star}(X),
w_{1},w_{2}\in \Phi(X)$ ($w_{1},w_{2}$ may be empty). By induction
we have $u'|_s=\sum_{i} \alpha_{i}u_{i}'|_{s_{i}}$ where each $
s_i\in S, u_{i}'|_{s_{i}}$ is normal $ s_{i}\mbox{-word}\
\mbox{and}\ \overline{u_{i}'|_{s_{i}}}\leq \overline{u'|_{s}}$. Thus
$u|_{s}=\sum_{i} \alpha_{i} u_{i}|_{s_{i}}$ where each
$u_{i}|_{s_{i}}=w_{1}P(u_{i}'|_{s_{i}})w_{2}$ is normal $s_{i}$-word
and
$u_{i}|_{\overline{s_{i}}}=w_1P(u_{i}'|_{\overline{s_{i}}})w_2\leq
\overline{w_1P(u'|_{s})w_2} =\overline{u|_{s}}. \hfill \blacksquare
$

\begin{lemma}\label{Deng6}
Let $S$ be a Gr\"{o}bner-Shirshov  basis  in  $RB(X), \ u_1, u_2\in
\Phi^\star(X)$
 and $s_1, s_2\in S$. If
 $w=u_1|_{\overline{s_1}}=u_2|_{\overline{s_2}}$
 where  each $u_i|_{s_i}$ is  a  normal $s_i$-word, $i=1,2$,  then
$$
u_1|_{s_1}\equiv u_2|_{s_2} \ \ mod (S,w).
$$
\end{lemma}
{\bf Proof:} \ There are three cases to consider. In fact, all cases
are essentially true in the same way.

(i)\ \  $\overline{s_1}$ and $\overline{s_2}$ are disjoint. In this
case, there exits a $(\star_1,\star_2)$-Rota-Baxter word $\Pi$ such
that
$$\Pi|_{\overline{s_1},\
\overline{s_2}}=u_1|_{\overline{s_1}}=u_2|_{\overline{s_2}}
$$
where $u_1|_{s_1}=\Pi|_{s_1, \overline{s_2}}$ is normal $s_1$-word
and $u_2|_{s_2}=\Pi|_{\overline{s_1}, s_2}$ is normal $s_2$-word.
Then
$$u_2|_{s_2}-u_1|_{s_1}
=\Pi|_{\overline{s_1}, \ s_2}- \Pi|_{s_1, \ \overline{s_2}}
=-\Pi|_{s_1-\overline{s_1},s_2} +\Pi|_{s_1,\ s_2-\overline{s_2}}.
$$
By Lemma \ref{lemma3.5}, we may suppose that
$$
-\Pi|_{s_1-\overline{s_1}, \ s_2}=\sum_{t}\alpha_{2t}u_{2t}|_{s_2},\
\Pi|_{s_1,s_2-\overline{s_2}}=\sum_{l}\alpha_{1l}u_{1l}|_{s_1},
$$
where $u_{2t}, u_{1l}\in \Phi^{\star}(X)$.

 Since  $S $ is a Gr\"{o}bner-Shirshov basis,
 by Lemma \ref{Deng5}, we have
$$
u_{2t}|_{s_2}=\sum_{n}\beta_{2tn}w_{2tn}|_{s_{_{2tn}}}
$$
and
$$
u_{1l}|_{s_1}=\sum_{m}\beta_{1lm}w_{1lm}|_{s_{_{1lm}}}
$$ where
$\beta_{1lm},\beta_{2tn}\in k, \ s_{_{2tn}},s_{_{1lm}}\in S,\
w_{2tn}|_{s_{_{2tn}}}\ \mbox{and}\ w_{1lm}|_{s_{_{1lm}}}$ are normal
$s_{_{2tn}}$- and $s_{_{1lm}}$-words, respectively,
$w_{2tn}|_{\overline{s_{_{2tn}}}}\leq
\overline{u_{2t}|_{s_{2}}}=\overline{\Pi|_{s_1-\overline{s_1},
s_2}}< \Pi|_{\overline{s_1},\overline{s_2}}=w$ and
$w_{_{1lm}}|_{\overline{s_{_{1lm}}}}<w.$ Therefore,
$$
u_2|_{s_2}-u_1|_{s_1}=\sum_{t,n}\alpha_{2t}\beta_{2tn}w_{2tn}|_{s_{_{2tn}}}+
\sum_{l,m}\alpha_{1l}\beta_{1lm}w_{1lm}|_{s_{_{1lm}}}.
$$
It
follows that
$$
u_1|_{s_1}\equiv u_2|_{s_2} \ \ mod (S,w).
$$

(ii)   $\overline{s_1}$ and $\overline{s_2}$ have nonempty
intersection but do not include each other. Without loss of
generality, we may assume that $\overline{s_1}a=b\overline{s_2}$ for
some $a, b\in\Phi(X)$.   This implies that $s_1a$ is normal
$s_1$-word and $bs_2$ is normal $s_2$-word. Then there exists a
$\Pi\in\Phi^\star(X)$ such that
$$
\Pi|_{\overline{s_1}a}=u_1|_{\overline{s_1}}=
u_2|_{\overline{s_2}}=\Pi|_{b\overline{s_2}}
$$
where $\Pi|_{s_1a}$ is normal $s_1a$-word and $\Pi|_{bs_2}$ is
normal $bs_2$-word. Thus, we have
$$
u_{_2}|_{s_2}-u_1|_{s_1} =\Pi|_{bs_2}-\Pi|_{s_1a}
=-\Pi|_{s_1a-bs_2}.
$$
Since  $S $ is a Gr\"{o}bner-Shirshov basis,
$$
s_1a-bs_2=\sum_{j}\alpha_jv_j|_{s_j}
$$
where each $\alpha_j\in k,  \ s_j\in S$,
$v_j|_{\overline{s_j}}<\overline{s_1}a$ and $v_j|_{s_j}$ is normal
$s_j$-word.  Let $ \Pi|_{v_j|_{s_j}}=\Pi_{j}|_{s_j} $. Then, by
Lemmas \ref{lemma3.5} and  \ref{Deng5},
$$\Pi_{j}|_{s_j}=\sum_{l}\beta_{jl}w_{jl}|_{s_{jl}}$$ where
each $\beta_{jl}\in k, \ s_{jl}\in S,\ w_{jl}|_{s_{jl}}$
is also normal $s_{jl}$-word, $w_{jl}|_{\overline{s_{jl}}}\leq
\overline{\Pi_j|_{s_j}}= \overline{\Pi|_{v_{j}|_{s_{j}}}}<
\Pi|_{\overline{s_1}a}=w$. Therefore
$$
u_2|_{s_2}-u_1|_{s_1}=\sum_{j}\alpha_j\Pi_{j}|_{s_j}=\sum_{j,l}\alpha_j\beta_{jl}w_{jl}|_{s_{jl}}.
$$
It follows that
$$
u_1|_{s_1}\equiv u_2|_{s_2} \ \ mod (S,w).
$$

(iii) One of  $\overline{s_1}$, $\overline{s_2}$ is contained in the
other. For example, let $ \overline{s_1} =u|_{\overline{s_2}}$ for
some $\star$-Rota-Baxter word $u$. Then
$$
w=u_2|_{\overline{s_2}}=u_1|_{u|_{\overline{s_2}}}
$$
and
$$
u_2|_{s_2}-u_1|_{s_1}=u_1|_{u|_{s_2}}-u_1|_{s_1}=-u_1|_{ s_1-u|_{s_2}}.
$$
Similar to  (ii), we can obtain the result. \hfill $\blacksquare$

\begin{lemma}\label{Deng7}
Let $S \subseteq RB(X)$ be a set of monic polynomials.
$Irr(S)\triangleq \{u\in \Phi(X)| u\neq v|_{\overline{s}}, s\in S,
v|_{s}\ $is normal $s$-word\}. Then for any $f\in RB(X),\ f$ has an
expression
$$
f=\sum_{i} \alpha_{i}u_{i}+\sum_{j} \beta_{j}v_{j}|_{s_{j}}
$$
where $\alpha_{i}, \beta_{j}\in k,\ u_{i}\in Irr(S),\
\overline{u_{i}}\leq \overline{f},\ s_j\in S, v_{j}|_{s_{j}}$ is
normal $s_j$-word and $v_{j}|_{\overline{s_{j}}}\leq \overline{f}$.
\end{lemma}
{\bf Proof:} Let $f=\sum_{i} \alpha_{i}u_{i},\ 0\neq \alpha_{i}\in
k, u_{1}> u_{2}> \cdots $. If $u_{1}\in Irr(S)$, then let
$f_{1}=f-\alpha_{1}u_{1}$. If $u_{1}\notin Irr(S)$, i.e., there
exists $ s_{1}\in S$ such that $u_{1}=v_{1}|_{\overline{s_{1}}}$
where $v_{1}|_{s_{1}}$ is  normal $s_{1}$-word, then let
$f_{1}=f-\alpha_{1}v_{1}|_{s_{1}}$. In both cases
$\overline{f_{1}}<\overline{f}$. Now, the result follows from
induction on $\overline{f}. \hfill \blacksquare$

\begin{theorem}(Composition-Diamond lemma for Rota-Baxter algebras)\label{Deng8}
Let $RB(X)$ be a free Rota-Baxter algebra over a field of
characteristic 0 and $S$  a set of monic polynomials in $RB(X)$, $
>$ the monomial ordering on $\Phi(X)$ defined as before and $Id(S)$ the Rota-Baxter ideal of
$RB(X)$ generated by $S$. Then the following statements are
equivalent.
\begin{enumerate}
\item[(I)] $S $ is a Gr\"{o}bner-Shirshov basis in $RB(X)$.
\item[(II)] $  f\in Id(S)\Rightarrow
\bar{f}=u|_{\overline{s}}$  for some $u \in \Phi^{\star}(X)$,\ $s\in
S$.
\item[(III)] $Irr(S) = \{ u\in \Phi(X) |  u \neq
v|_{\bar{s}}, s\in S, v|_{s}\ \mbox{is normal}\  s\mbox{-word}\}$ is
a $k$-basis of $RB(X|S)\\=RB(X)/Id(S)$.
\end{enumerate}
\end{theorem}

{\bf Proof:}\ \ (I)$\Longrightarrow$ (II)\ \ Let  $0\neq f\in
Id(S)$. By Lemmas \ref{lemma3.5} and \ref{Deng5}, we may assume that
$f=\sum\limits_{i=1}^{n}\alpha_i u_i|_{s_i}$ where each $\alpha_i\in
k$, $s_i\in S$ and $u_i|_{s_i}$ is normal $s_i$-word. Let
$w_i=u_i|_{\overline{ s_i}}$. We arrange these leading words in
non-increasing order by
$$
w_1= w_2=\cdots=w_m >w_{m+1}\geq \cdots\geq w_n.
$$
Now we prove the result by  induction  on $m$ and $w_1$.

If $m=1$, then $\bar{f}=u_1|_{\overline{s_1}}$ and there is nothing
to prove. Now we assume that $m\geq 2$. Then $ u_1|_{\overline{
s_1}}=w_1=w_2=u_2|_{\overline{s_2}}$. Since $S $ is a
Gr\"{o}bner-Shirshov basis  in $RB(X)$, by Lemma \ref{Deng6}, we
have
$$
u_2|_{s_2}-u_1|_{s_1}=\sum_{j}\beta_jv_j|_{s_j}
$$
where each $\beta_j\in k, \ s_j\in S, \ v_j|_{\overline{s_j}}<w_1$
and $v_j|_{s_j}$ is normal $s_j$-word. Therefore, since
$$
\alpha_1u_1|_{s_1}+\alpha_2u_2|_{s_2}
=(\alpha_1+\alpha_2)u_1|_{s_1}+\alpha_2(u_2|_{ s_2} -u_1|_{ s_1}),
$$
we have
$$
f=(\alpha_1+\alpha_2)u_1|_{s_1}+\sum_{j}\alpha_2\beta_jv_j|_{s_j}+
\sum\limits_{i=3}^{n}\alpha_iu_i|_{ s_i}.
$$

If either $m>2$ or $\alpha_1+\alpha_2\neq 0$, then the result
follows from induction on $m$. If $m=2$ and $\alpha_1+\alpha_2=0$,
then the result follows from induction on $w_1$.

(II)$\Longrightarrow$ (III) For any $f\in RB(X)$, by Lemma
\ref{Deng7},  $f+Id(S)$ can be expressed by the elements of
$Irr(S)+Id(S)$. Now Suppose that
$\alpha_1u_1+\alpha_2u_2+\cdots+\alpha_nu_n=0$ in $ RB(X|S)$ with
each $u_i\in Irr(S)$,\  $u_1>u_2>\cdots>u_n$\  and each $0\neq
\alpha_i\in k$. Then, in $RB(X)$,
$$
g=\alpha_1u_1+\alpha_2u_2+\cdots+\alpha_nu_n\in Id(S).
$$
By (II), we have $u_1=\bar{g}\notin Irr(S)$, a contradiction. Hence,
$Irr(S)$ is $k$-linearly independent. This shows that $Irr(S)$ is a
$k$-basis of $ RB(X|S)$.

(III)$\Longrightarrow $(I) For  any $(f,g)_w$ where $f,g \in S$, by
Lemma \ref{Deng7}, $(f,g)_w=\sum_{i} \alpha_{i}u_{i}+\sum_{j}
\beta_{j}v_{j}|_{s_{j}}$ where  $\alpha_{i}, \beta_{j}\in k,
u_{i}\in Irr(S), s_j\in S, \overline{u_{i}}\leq
\overline{(f,g)_w}<w, v_{j}|_{s_{j}}$ is normal $s_{j}$-word and
$v_{j}|_{\overline{s_{j}}}\leq \overline{(f,g)_w}<w$. Since
$(f,g)_w\in Id(S)$ and by (III), each $\alpha_{i}=0$. Thus
$(f,g)_w\equiv0\ mod(S,w)$.

For any composition of left multiplication, $P(u)\cdot f$ where
$f\in S$,  by Lemma \ref{Deng7} and by (III), we have $P(u)\cdot f =
\sum_{j} \beta_{j}v_{j}|_{s_{j}}$ where each $\beta_{j}\in k, \
s_j\in S, v_{j}|_{s_{j}}$ is normal $s_{j}$-word and
$v_{j}|_{\overline{s_{j}}}\leq \overline{P(u)\cdot f}$. This shows
that each composition of left multiplication in $S$ is trivial
$mod(S)$. Similarly, each composition of right multiplication in $S$
is trivial $mod(S)$. $\hfill \blacksquare$

\section{Applications}

In \cite{R069, Ca}, G. Rota and P. Cartier constructed free
commutative Rota-Baxter algebra with certain restriction. While in
\cite{GKg00a, GKg00b}, L. Guo and W. Keigher constructed free
commutative Rota-Baxter algebra with unity by mixable shuffle
algebra, also see \cite{GX}. In this section, by using the
Composition-Diamond lemma for Rota-Baxter algebras (Theorem
\ref{Deng8}), we have the following results: 1. We give a
Gr\"{o}bner-Shirshov basis for the free commutative Rota-Baxter
algebra $RB(X|S)$ and as an application, a linear basis of the free
commutative Rota-Baxter algebra without unity is obtained. 2. We
show that every countably generated Rota-Baxter algebra with weight
0 can be embedded into a two-generated Rota-Baxter algebra. 3. We
prove the $\frac{1}{2}$-PBW Theorems for dendriform dialgebra and
dendriform trialgebra.

\begin{theorem}\label{Deng9} Let $I$ be a well ordered set,
$X=\{x_i|i\in I\}$ and the ordering on $\Phi(X)$ defined as before.
Let
\begin{eqnarray}
\label{e1}&&f=x_{i}x_{j}-x_{j}x_{i}, \ \ i>j, \ i,j\in I\\
\label{e2}&&g=P(u)x_{i}-x_{i}P(u), \ \  u\in \Phi(X),\ i\in I
\end{eqnarray}
Let $S$ consist of (\ref{e1}) and (\ref{e2}). Then $S$ is a
Gr\"{o}bner-Shirshov basis in $RB(X)$.
\end{theorem}
{\bf Proof:} Denote by $s\wedge l$ the composition of the
polynomials of type $s$ and type $l$.

The ambiguities of all possible compositions of the polynomials  in
$S$
 are only as below:
 \begin{description}
  \item $1\wedge1\ \ \ x_{i}x_{j}x_{l}, \ \ i>j>l;$
  \item $2\wedge1\ \ \  P(u)x_{i}x_{j}\  \mbox{and} \ \ P(u)x_{l}=
  P(v|_{x_{i}x_{j}})x_{l}, \ \ u\in \Phi(X),\ v\in \Phi^{\star}(X), i>j;$
  \item $2\wedge2\ \ \ P(u)x_{i}=P(v|_{P(w)x_{j}})x_{i},\ \
  u,w\in \Phi(X), v\in \Phi^{\star}(X).$
\end{description}

All possible compositions of left multiplication are $P(v)\cdot g$
where $v\in \Phi(X)$.

Now we check that all the compositions are trivial. We only check
$1\wedge1$ as an example.

For $1\wedge1$, let $f=x_{i}x_{j}-x_{j}x_{i},
g=x_{j}x_{l}-x_{l}x_{j}$. Then $w=x_{i}x_{j}x_{l}$.
\begin{eqnarray*}
(f,g)_w&=& (x_{i}x_{j}-x_{j}x_{i})x_{l}-x_{i}(x_{j}x_{l}-x_{l}x_{j})\\
&=& x_{i}x_{l}x_{j}-x_{j}x_{i}x_{l}\\
&\equiv& x_{l}x_{j}x_{i}-x_{l}x_{j}x_{i}\\
&\equiv& 0 \ \ \ mod(S, x_{i}x_{j}x_{l}).
\end{eqnarray*}

For any $P(v)$, let  $ P(v)=P^{l}(v')$ where $l\geq 1, v'\notin
P(\Phi(X))$.  Then
\begin{eqnarray*}
&&P(v)(P(u)x_{i}-x_{i}P(u))\\
&=&P^{l}(v')\cdot (P^{t}(u')x_{i}-x_{i}P^{t}(u'))\\
&\equiv& (P^{l}(v')\cdot P^{t}(u'))x_{i}-x_{i}(P^{l}(v')\cdot
P^{t}(u'))\ \ mod(S)
\end{eqnarray*}
where $l,t\geq 1, u',v'\notin P(\Phi(X))$. By Lemmas \ref{Deng1} and
\ref{Deng2},
$$
P^{l}(v')\cdot P^{t}(u')=\alpha_{(l,t,1)}P^{l+t-1}(P(v')\cdot u')+
\varepsilon \triangleq \alpha_0\varepsilon_{0}+ \sum_{j=1}
\alpha_{j}\varepsilon_{j}
$$
where $\varepsilon\in RB(X)$, each $\varepsilon_{j}\in P(\Phi(X)),\
\varepsilon_{j}<\overline{P^{l+t-1}(P(v')\cdot
u')}=\overline{P^{l}(v')\cdot P^{t}(u')}, \alpha_j\in k,
j=1,2,\ldots , \ \varepsilon_{0}=\overline{P^{l+t-1}(P(v')\cdot
u')}=\overline{P^{l}(v')\cdot P^{t}(u')}\in P(\Phi(X))\  \mbox{and}\
\alpha_0=\alpha_{(l,t,1)}.$ Therefore
$$
P(v)(P(u)x_{i}-x_{i}P(u))\equiv\alpha_0(\varepsilon_{0}x_{i}-x_{i}\varepsilon_{0})
+\sum_{j=1}\alpha_{j}(\varepsilon_{j}x_{i}-x_{i}\varepsilon_{j}) \ \
\  mod(S)$$ where each $\varepsilon_{m}x_{i}-x_{i}\varepsilon_{m}\in
S$ and $\overline{\varepsilon_{m}x_{i}-x_{i}\varepsilon_{m}}\leq
\overline{P^{l}(v')\cdot (P^t(u')x_i)}=\overline{P(v)\cdot
(P(u)x_i)}, m\geq 0$. Thereby
$$
P(v)(P(u)x_{i}-x_{i}P(u))\equiv 0\ \ \ mod(S).
$$
Thus all the compositions in $S$ are trivial.

This shows that $S$ is a Gr\"{o}bner-Shirshov basis. $ \hfill \
\blacksquare$

 \ \

We now use the same notations as in Theorem \ref{Deng9}. Clearly,
$RB(X|S)$ is  the free commutative Rota-Baxter algebra generated by
$X$ with weight $\lambda$.

The following corollary follows from our Theorems \ref{Deng8} and
\ref{Deng9}.

\begin{corollary}\label{Deng10}(\cite{GKg00b}) Let $RB(X|S)$ be the free
commutative Rota-Baxter algebra generated by $X$ with weight
$\lambda$ as in Theorem \ref{Deng9}. Then
$$
Irr(S)=Y_1\cup Y_2\cup Y_3$$
 is a $k$-basis of   $RB(X|S)$ where
\begin{eqnarray*}
 Y_1 &=& \{x_{1}x_{2}\cdots x_{n}\in
S(X)|n\in {\cal N^+}, x_i\in X,
x_{1}\leq x_{2}\leq \cdots \leq x_{n}\},\\
Y_2 &=& \{P^{l_{1}}(u_{1}P^{l_{2}}(\cdots (u_{t-1}P^{l_{t}}(u_{t})
)\cdots ))\in
\Phi(X)|u_j\in Y_1, t\geq1, l_j\geq1, j=1,2,\dots,t\},\\
Y_3 &=&\{uv\in \Phi(X)|u\in Y_1,v\in Y_2\}.
\end{eqnarray*}$\hfill \ \blacksquare$
\end{corollary}

Now,  we show that every countably generated Rota-Baxter algebra
with weight 0 can be embedded into a two-generated Rota-Baxter
algebra.

\begin{theorem}\label{Deng12}
Every countably generated Rota-Baxter algebra with weight 0 can be
embedded into a two-generated Rota-Baxter algebra.
\end{theorem}

\textbf{Proof:}\  Since every countably generated Rota-Baxter
algebra $A$ is countably dimensional, we may assume that $A$ has a
$k$-basis $X=\{x_{i}| i=1,2,\ldots\}$. Then $A$ can be expressed as
$A=RB(X|S)$, where $
S=\{x_{i}x_{j}=\{x_{i},x_{j}\},P(x_{i})=\{P(x_{i})\}\mid i,j\in
{\cal N^+}\}$ and for any $u\in \Phi(X)$, $\{u\}$ means the linear
combination of $u$ by $x_t, \ x_t\in X$.

For $u\in\Phi(X)$, we denote $supp(u)$  the set of $x_i\in X$
appearing in the word $u$, for example, if $u=x_1P(x_2)x_3x_1$, then
$supp(u)=\{x_1, x_2, x_3\}$.

Let $H=RB(X,a,b| S_{1})$, where $S_{1}$ consists of the following
relations:
\begin{enumerate}
\item[I.]\ $x_{i}x_{j}=\{x_{i}x_{j}\},\ i,j\in {\cal N^+},$
\item[II.]\ $P(x_{i})=\{P(x_{i})\},\ i\in {\cal N^+},$
\item[III.]\ $aab^{i}ab=x_{i},\ i\in {\cal N^+},$
\item[IV.]\ $P(t)=0,\ t\in \Phi'(X,a,b)$,
\end{enumerate}
where $\Phi'(X,a,b)=\{u\in\Phi(X,a,b) | \exists u'\in\Phi(X,a,b),
u'\in Irr(\{III\}),\ supp(u')\cap\{a,b\}\neq\emptyset, s.t. ,\
u'\equiv u\ mod(\{III\},u)\}$.

Let the ordering  on $\Phi(X,a,b)$ be defined as before, where
$x_i>a>b, \ i\in {\cal N^+}$.

We will prove that $S_{1}$ is a Gr\"{o}bner-Shirshov basis in
$RB(X,a,b)$.

The ambiguities $w$ of all possible compositions of in $S_{1}$ are:
$$1)\ \ x_{i}x_{j}x_{l}\ \   2)\ \ P(t|_{x_{i}x_{j}}) \ \ 3)\ \
P(t|_{P(x_{i})}) \ \ 4)\ \ P(t|_{P(t')}) \ \ 5)\ \
P(t|_{aab^{i}ab})$$ where $x_{i},x_{j},x_{l}\in X, t, t'\in
\Phi'(X,a,b)$.

All possible compositions of left and right multiplications in $S_1$
are as below:
\begin{enumerate}
\item[6)]\ \ $P(u)(P(x_i)-\{P(x_i)\}), i\in {\cal N^+}, u\in \Phi(X,a,b),$
\item[7)]\ \ $(P(x_i)-\{P(x_i)\})P(u), i\in {\cal N^+}, u\in \Phi(X,a,b),$
\item[8)]\ \ $P(u)(P(t)-0),u\in \Phi(X,a,b), t\in  \Phi'(X,a,b),$
\item[9)]\ \ $(P(t)-0)P(u), u\in \Phi(X,a,b), t\in  \Phi'(X,a,b).$
\end{enumerate}

We have to check that all these compositions are trivial. We just
take 5) and 6) as examples.

For 5), let $f=P(t|_{aab^{i}ab}),\  g=aab^{i}ab-x_{i},\  t
 \in \Phi'(X,a,b),\ x_{i}\in X $. Then $w=P(t|_{aab^{i}ab})$ and
$$
(f,g)_{w}=P(t|_{x_{i}})\equiv0 \ mod(\{IV\},w).
$$

For 6),
\begin{eqnarray*}
&&P(u)(P(x_i)-\{P(x_i)\})\\
&=& P(P(u)x_i)+P(u P(x_i))-P(u)\{P(x_i)\}\\
&\equiv&P(P(v)x_i)+P(v P(x_i))-P(v)\{P(x_i)\} \ \ \ mod(\{III\})
\end{eqnarray*}
where $v\in\Phi(X,a,b)$ such that $u\equiv v \ mod(\{III\},u)$.
There are two cases to consider.
\begin{description}
\item{a):} If $v\in \Phi'(X,a,b)$, then $mod(\{IV\}), \ P(P(v)x_i)\equiv0,
\ P(v P(x_i))\equiv0$ and $P(v)\{P(x_i)\}\equiv0$. Therefore,
$$
P(u)(P(x_i)-\{P(x_i)\})\equiv 0\ \ mod(S_1).
$$
\item{b):} If $v\not\in \Phi'(X,a,b)$, then $v\in \Phi(X)$.
Since $v\equiv \sum\beta_jx_j \ \ mod(\{I,II\})$,
\begin{eqnarray*}
P(u)(P(x_i)-\{P(x_i)\})\equiv 0\ \  mod(S_1).
\end{eqnarray*}
\end{description}

Thus, $S_{1}$ is a Gr\"{o}bner-Shirshov basis in $ RB(X,a,b)$. By
Theroem \ref{Deng8}, $A$ can be embedded into the Rota-Baxter
algebra $H$ which is generated by $\{a,b\}$. \hfill $\blacksquare$

\ \

Recall that a dendriform dialgebra (see \cite{Lo}) is a $k$-module
$D$ with two binary operations $\prec$ and $\succ$ such that
\begin{enumerate}
\item \ \ \  $(x\prec y)\prec z=x\prec(y\prec z+y\succ z)$
\item \ \ \ $(x\succ y)\prec z=x\succ(y\prec z)$
\item \ \ \  $(x\prec y+x\succ y)\succ z=x\succ(y\succ z)$
\end{enumerate}
for any $x,y,z\in D.$

A dendriform trialgebra (see \cite{LR}) is a $k$-module $T$ equipped
with three binary operations $\prec, \succ$ and $\circ$ that satisfy
the relations
\begin{enumerate}
\item\ \ \ $(x\prec y)\prec z=x\prec(y*z)$
\item\ \ \ $(x\succ y)\prec z=x\succ(y\prec z)$
\item\ \ \ $(x*y)\succ z=x\succ(y\succ z)$
\item\ \ \ $(x\succ y)\circ z=x\succ(y\circ z)$
\item\ \ \ $(x\prec y)\circ z=x\circ(y\succ z)$
\item\ \ \ $(x\circ y)\prec z=x\circ(y\prec z)$
\item\ \ \ $(x\circ y)\circ z=x\circ (y\circ z)$
\end{enumerate}
for any $x,y,z\in T$ where $x*y=x\prec y+x\succ y+x\circ y.$

Denote by $D(X|S)\ (T(X|S))$ the dendriform dialgebra (trialgebra)
generated by $X$ with defining relations $S$.

Suppose that $(D, \prec,\succ)$ is a dendriform dialgebra with a
linear basis $X=\{x_n|n\in I\}$. Let $x_i\prec x_j=[x_i,x_j],
x_i\succ x_j=[[x_i,x_j]]$, where $[x_i,x_j]$ and $[[x_i,x_j]]$ are
linear combinations of $x_n\in X$.

Then $D$ has the expression
$$
D=D(X|x_i\prec x_j=[x_i,x_j], x_i\succ x_j=[[x_i,x_j]]).
$$
In \cite{EG08a}, K. Ebrahimi-Fard and L. Guo proved that
$$
U(D)=RB(X|x_iP(x_j)+\lambda x_ix_j=[x_i,x_j], P(x_i)x_j=[[x_i,x_j]])
$$
is the universal enveloping algebra of $D$, where $\lambda \in k$.

\begin{theorem}\label{Deng13} Suppose that $(D, \prec,\succ)$ is a dendriform dialgebra
with a linear basis $X=\{x_n|n\in I\}$, where $I$ is a well ordered
set. Let $x_i\prec x_j=[x_i,x_j], x_i\succ x_j=[[x_i,x_j]]$, where
$[x_i,x_j]$ and $[[x_i,x_j]]$ are linear combinations of $x_n\in X$.
Let $\lambda\in k, \lambda\neq0$. Let $S$ be the set consisting of
\begin{eqnarray}
\label{e3}&&f_1=x_iP(x_j)+\lambda x_ix_j-[x_i,x_j],\\
\label{e4}&&f_2=P(x_i)x_j-[[x_i,x_j]],\\
\label{e5}&&f_3=x_ix_jx_l-{\lambda}^{-1}[x_i,x_j]x_l+{\lambda}^{-1}x_i[[x_j,x_l]].
\end{eqnarray}
Then in $RB(X)$, all intersection compositions and inclusion
compositions in $S$
 are trivial.
\end{theorem}
{\bf Proof:} The ambiguities of all possible compositions of the
polynomials in $S$ are only as below:
 \begin{description}
  \item $3\wedge4\ \ \  w_{34}=x_{i}P(x_{j})x_{l};\ \ \ 4\wedge3\ \ \  w_{43}=P(x_{i})x_{j}P(x_l);$
  \item $4\wedge5\ \ \ w_{45}=P(x_{l})x_ix_jx_s;\ \ \ \ 5\wedge3\ \ \ w_{53}=x_ix_jx_sP(x_l);$
  \item $5\wedge5 \ \ \ w_{55}=x_ix_jx_sx_l \ \ \ \ \ \mbox{and}\ \ \ \ \ w'_{55}=x_ix_jx_hx_lx_s.$
\end{description}
We only check $4\wedge3$ as example. The other cases are similar.
Let $f=P(x_{i})x_{j}-[[x_i,x_{j}]], g=x_{j}P(x_{l})+\lambda
x_{j}x_{l}-[x_j,x_l]$. Then $w_{43}=P(x_{i})x_{j}P(x_{l})$ and
\begin{eqnarray*}
(f,g)_{w_{43}}&=& (P(x_{i})x_{j}-[[x_i,x_{j}]])P(x_{l})-P(x_{i})(x_{j}P(x_{l})+\lambda x_{j}x_{l}-[x_j,x_l])\\
&=& -\lambda P(x_{i})x_{j}x_{l}+P(x_{i})[x_{j},x_{l}]-[[x_i,x_j]]P(x_l)\\
&\equiv& -\lambda [[x_{i},x_{j}]]x_{l}+[[x_{i},[[x_{j},x_{l}]]]]-(-\lambda[[x_i,x_j]]x_l+[[[x_i,x_j]],x_l])\\
&\equiv& -[[[x_i,x_j]],x_l]+[[x_{i},[[x_{j},x_{l}]]]]\\
&\equiv&0\ \ \ mod(S, P(x_{i})x_{j}P(x_l))
\end{eqnarray*}
by using the relation
$
(x\succ y)\prec z=x\succ(y\prec z).
$

Thus, the theorem is proved. $\hfill \blacksquare$

\ \

Suppose that $(T, \prec,\succ,\circ)$ is a dendriform trialgebra
with a linear basis $X=\{x_n|n\in I\}$. Let $x_i\prec x_j=[x_i,x_j],
x_i\succ x_j=[[x_i,x_j]],x_i\circ x_j=\langle x_i,x_j\rangle$, where
$[x_i,x_j],[[x_i,x_j]]$ and $\langle x_i,x_j\rangle$ are linear
combinations of $x_n\in X$.

Then $T$ has the expression
$$
T=T(X|x_i\prec x_j=[x_i,x_j], x\succ x_j=[[x_i,x_j]], x_i\circ
x_j=\langle x_i,x_j\rangle).
$$

We notice in \cite{EG08a}, K. Ebrahimi-Fard and L. Guo have proved
that
$$
U(T)= RB(X|x_iP(x_j)=[x_i,x_j],P(x_i)x_j=[[x_i,x_j]],
 {\lambda}x_ix_j=\langle x_i,x_j\rangle)
$$
is the universal enveloping algebra of $T$, where $\lambda \in k$ .

\begin{theorem}\label{Deng14} Suppose that $(T, \prec,\succ,\circ)$ is a dendriform
trialgebra with a linear basis $X=\{x_n|n\in I\}$, where $I$ is a
well ordered set. Let $x_i\prec x_j=[x_i,x_j], x_i\succ
x_j=[[x_i,x_j]],x_i\circ x_j=\langle x_i,x_j\rangle$, where
$[x_i,x_j],[[x_i,x_j]]$ and $\langle x_i,x_j\rangle$ are linear
combinations of $x_n\in X$. Let $S$ be the set consisting of
\begin{eqnarray}
\label{e6}&&f_4=x_iP(x_j)-[x_i,x_j],\\
\label{e7}&&f_5=P(x_i)x_j-[[x_i,x_j]],\\
\label{e8}&&f_6={\lambda}x_ix_j-\langle x_i,x_j\rangle.
\end{eqnarray}
Then in $RB(X)$, all intersection compositions and inclusion
compositions in $S$
 are trivial.
\end{theorem}
{\bf Proof:} The ambiguities of all possible compositions of the
polynomials in $S$ are only as below:
 \begin{description}
  \item $6\wedge7\ \ \  w_{67}=x_{i}P(x_{j})x_{l};\ \ \ 7\wedge6\ \ \  w_{76}=P(x_{i})x_{j}P(x_l);$
  \item $7\wedge8\ \ \ w_{78}=P(x_{l})x_ix_j;\ \ \ \ 8\wedge6\ \ \ w_{86}=x_ix_jP(x_l);$
  \item $8\wedge8 \ \ \ w_{88}=\lambda x_ix_jx_l.$
\end{description}
We only check $7\wedge6$ as example. The other cases are similar.

Let $f=P(x_{i})x_{j}-[[x_i,x_{j}]], g=x_{j}P(x_{l})-[x_j,x_l]$. Then
$w_{76}=P(x_{i})x_{j}P(x_{l})$ and
\begin{eqnarray*}
(f,g)_{w_{76}}&=& (P(x_{i})x_{j}-[[x_i,x_{j}]])P(x_{l})-P(x_{i})(x_{j}P(x_{l})-[x_j,x_l])\\
&=& P(x_{i})[x_{j},x_{l}]-[[x_i,x_j]]P(x_l)\\
&\equiv& [[x_{i},[x_{j},x_{l}]]]-[[[x_i,x_j]],x_l]\\
&\equiv&0 \ \ \ mod(S, P(x_{i})x_{j}P(x_l))
\end{eqnarray*}
 by using the relation
$ (x\succ y)\prec z=x\succ(y\prec z). $

Hence, the proof is completed. $\hfill \blacksquare$

\ \

\end{document}